\documentclass[10pt]{article}
\usepackage{upgreek}
\usepackage[russian,greek,english]{babel}
\usepackage{amsmath}
\usepackage{amssymb}
\usepackage{amsthm}
\usepackage{amsfonts}
\usepackage{graphicx}
\usepackage{wrapfig}
\usepackage{mathtools}
\usepackage{tikz}
\usetikzlibrary{decorations.markings,fadings}
\usepackage{txfonts}
\usepackage{mathrsfs}
\usepackage{enumitem}
\usepackage[labelsep=period]{caption}
\usepackage{setspace}
\usepackage[pagewise]{lineno}
\usepackage{fancyhdr}

\tikzset{->-/.style={decoration={
  markings,
  mark=at position #1 with {\arrow{>}}},postaction={decorate}}}

\newcommand{\intl}[2]{\!\underset{#1}{\overset{#2}{\rotatebox[origin=rc]{15}{\large\ensuremath{\int}}}}}
\newcommand{\intL}[2]{\!\!\underset{#1}{\overset{#2}{\rotatebox[origin=rc]{15}{\Large\ensuremath{\int}}}}\!}
\newcommand{\inth}[2]{\!\!\!\underset{#1}{\overset{#2}{\rotatebox[origin=rc]{15}{\huge\ensuremath{\int}}}}\!\!}

\newcommand{\z}{\textsl{z}}

\newcommand{\D}{\partial}

\newcommand{\deq}{\overset{\operatorname{def}}{=}}

\renewcommand{\alpha}{\alphaup}

\renewcommand{\beta}{\betaup}
\renewcommand{\gamma}{\gammaup}
\renewcommand{\delta}{\deltaup}

\renewcommand{\epsilon}{\varepsilonup}
\renewcommand{\varepsilon}{\epsilonup}

\renewcommand{\zeta}{\zetaup}

\renewcommand{\eta}{\etaup}
\renewcommand{\theta}{\thetaup}
\renewcommand{\vartheta}{\varthetaup}

\renewcommand{\iota}{\iotaup}

\renewcommand{\kappa}{\varkappa}
\renewcommand{\lambda}{\lambdaup}

\renewcommand{\mu}{\muup}

\renewcommand{\nu}{\nuup}
\renewcommand{\xi}{\xiup}

\renewcommand{\pi}{\piup}

\renewcommand{\rho}{\rhoup}
\renewcommand{\varrho}{\varrhoup}
\renewcommand{\sigma}{\sigmaup}
\renewcommand{\varsigma}{\varsigmaup}

\renewcommand{\tau}{\tauup}
\renewcommand{\Upsilon}{\textrm{\greektext U}}
\renewcommand{\upsilon}{\upsilonup}
\renewcommand{\phi}{\upvarphi}
\renewcommand{\varphi}{\phiup}

\renewcommand{\chi}{\chiup}
\renewcommand{\psi}{\textrm{\greektext y}}
\renewcommand{\omega}{\omegaup}

\renewcommand{\mathbb}{\varmathbb}

\newtheorem*{definition}{Definition}
\newtheorem{theorem}{Theorem}
\newtheorem{lemma}{Lemma}

\newcommand\blfootnote[1]{
  \begingroup
  \renewcommand\thefootnote{}\footnote{#1}
  \addtocounter{footnote}{-1}
  \endgroup
}

\title{Exact travelling solution for a reaction-diffusion system with a piecewise constant production supported by a codimension-1 subspace}
\author{Anton S. Zadorin}
\date{\small\itshape
	Chimie Biologie Innovation, ESPCI Paris, CNRS, PSL University, 75005 Paris, France.\\
	Center for Interdisciplinary Research in Biology (CIRB), Coll\`ege de France, CNRS, INSERM, PSL
	Research University, Paris, France.\\
	Current affiliation: Max Planck Institute for Mathematics in the Sciences, Leipzig, Germany.}

\begin{document}
\maketitle

\begin{abstract} A generalisation of reaction diffusion systems and their travelling solutions to cases when the productive part of the reaction 
happens only on a surface in space or on a line on plane but the degradation and the diffusion happen in bulk are important for modelling various 
biological processes. These include problems of invasive species propagation along boundaries of ecozones, problems of gene spread in such situations, 
morphogenesis in cavities, intracellular reaction etc. Piecewise linear approximations of reaction terms in reaction-diffusion systems often result in 
exact solutions of propagation front problems. This article presents an exact travelling solution for a reaction-diffusion system with a piecewise 
constant production restricted to a codimension-1 subset. The solution is monotone, propagates with the unique constant velocity, and connects the 
trivial solution to a nontrivial nonhomogeneous stationary solution of the problem. The properties of the solution closely parallel the properties of 
monotone travelling solutions in classical bistable reaction-diffusion systems. \end{abstract}


\section*{Introduction}

\subsection*{Notations}

The following notational convention is used in the article.

\begin{itemize}

\item
$\mathbb R^n(x_1,\ldots,x_n)$: The space $\mathbb R^n$ with its usual topology coordinatised with variables $x_1$, \ldots, $x_n$.

\item
$\mathcal S(\mathbb R^n)$: The Schwartz's space (the complex-valued rapidly decreasing test functions). Abbreviated as $\mathcal S$ if obvious.

\item
$\mathcal S'(\mathbb R^n)$: The space of complex-valued tempered dsitributions on $\mathbb R^n$ (the dual of $\mathcal S(\mathbb R^n)$). Abbreviated as $\mathcal S'$ if obvious.

\item
$\mathcal O_M(\mathbb R^n)$: The multipliers of $\mathcal S(\mathbb R^n)$ and $\mathcal S'(\mathbb R^n)$, the space of slowly increasing functions\\ $\{f \in C^\infty(\mathbb R^n): \forall \sigma \in \mathbb N^n\, \exists k \in \mathbb N\, (1 + |x|^2)^{-k} |\D^\sigma f(x)| \to 0, |x| \to \infty\}$. Abbreviated as $\mathcal O_M$ if obvious.

\item
$\Omega$: $\mathbb R^2(x,y)$.

\item
$L$: $\mathbb R(x)$ understood as the special line $y = 0$ in $\Omega$, the usual inclusion $L \subset \Omega$ is implied.

\item
$\Omega_t$: $\mathbb R^3(t,x,y) = \mathbb R(t) \times \Omega$.

\item
$L_t$: $\mathbb R^2(t,x) = R(t) \times L$, the usual inclusion $L_t \subset \Omega_t$ is implied.

\item
$f*g$: The convolution of $f,g \in \mathcal S'$.

\item
$K_0$: the Macdonald function.

\end{itemize}

\subsection*{Biological motivation}

Many systems in different branches of biology, primarily in cellular biology and in spatial population dynamics, can be idealised by nonlinear chemical 
reactions supported by codimension-1 subsets coupled to free transport (typically by diffusion) and maybe linear degradation in bulk. Direct examples 
would be biochemical reactions catalysed by membrane-bound enzymes in cells or in organelles, chemical interaction of different parts of an epithelium 
separated by a lumen, chemical communication of bacteria in biofilms etc. For example in morphogenesis, many embryological processes that rely on 
diffusing morphogens happen in a tissue that surrounds or is surrounded by some sort of cavity or lumen: gastrulation, early neurulation, formation of 
cerebral vesicles, angiogenesis to name a few.

At a different scale, under some circumstances population dynamics can be approximated as a reaction-diffusion process. Historically, the discovery 
of travelling front solutions of nonlinear diffusion equations came from considering a problem of spatial spreading of an advantageous gene in a 
population of organisms \cite{Fisher1937,Kolmogorov1937}. In case of a not very mobile short living species that inhabits an edge between two ecozones 
(an edge of a forest, a coast etc.) that spreads dispersing propagula (larvae, seeds, spores) that are carried everywhere but can develop/sprout 
only on the edge, the dynamics of spreading of such a species can be modelled by a reaction-diffusion system with the nonlinear reaction supported 
by a line (the edge). It is known that many invasive species establish first along edges like river banks, roadsides, and other disturbed landscapes.  
The same mathematical structure describes the propagation of a new advantageous gene in a population of such organisms that is already well established 
on the edge. In that case, the organisms themselves do not have to spread via mobile propagula, but the long range dispersion is carried out by their 
gametes/pollen. An archetypal example is a littoral animal with low mobility that relies on a planktonic larva for spreading: shallow water 
barnacles, gastropods, bivalves, decapods, corals etc.

Despite the widespreadness of nonlinear reactions supported by codimension-1 subsets their theoretical understanding is limited. Usually the 
understanding is based on an analogy with pure reaction-diffusion systems, where the coupling of the reaction rates on the reaction supporting subset 
through the volume diffusion is approximated in an \emph{ad hoc} manner either by averaging or by integration over the volume in the direction 
transverse to the subset. Such approaches, however, neglect the nonlocal character of the problem that is a consequence of the diffusive coupling 
through the bulk.

In a recent series of works by H. Berestycki, J.M. Roquejoffre, and L. Rossi \cite{Berestycki2013a,Berestycki2013b,Berestycki2016}, a closely related 
problem is considered. They studied a general behaviour of a travelling solution to a reaction-diffusion problem that describes an invasion of a species 
in a plane facilitated by a fast diffusion on a line.

In this article, I present an exactly solvable model of the propagation problem in a reaction-diffusion system where the growth happens only on a 
plane (in space) or on a line (on plane) and the growth rate is approximated by a piecewise constant function.

\subsection*{Informal problem statement}

\begin{figure}[!t]
	\begin{tikzpicture}
		\begin{scope}
			\node at (-1.5,2) {\textbf{A}};
			\draw[dotted] (0,-1) -- (2,-2) -- (4,-1) -- (2,0) -- (0,-1);
			\draw[fill=gray,loop] (0,0) -- (1.5,-0.75) -- (3.5,0.25) -- (2,1);
			\draw[fill=lightgray,loop] (1.5,-0.75) -- (2,-1) -- (4,0) -- (3.5,0.25);
			\draw[thick,->] (0,0) -- (2,-1) node[below left] {$x$};
			\draw[thick,->] (0,0) -- (2,1) node[above left] {$\z$};
			\draw[thick,->] (0,-1) -- (0,1) node[left] {$y$};
			\draw[dotted] (0,1) -- (2,0) -- (4,1) -- (2,2) -- (0,1);
			\draw[dotted] (2,-2) -- (2,0);
			\draw[dotted] (4,-1) -- (4,1);
			\draw[dotted] (2,1) -- (2,2);
			\node [rotate=26.6,white] at (2.2,-0.1) {\footnotesize\bfseries activity front};
			\draw[thick,->] (0.1,-0.3)node[left]{\parbox{2cm}{\centering\footnotesize propagation of activity with a constant speed}}
				-- node[below left]{$\varv$} (1,-0.75);
			\draw[->] (4,-0.5) node[right] {\parbox{2cm}{\centering\footnotesize
				nonlinear growth on a plane}} -- (3.2,-0.5);
			\draw[<-] (3.3,1.6) -- (4,1.8) node[right] {
				\parbox{2.5cm}{\centering\footnotesize
				diffusion and linear degradation in bulk}};
		\end{scope}

		\begin{scope}[xshift=9cm]
			\node at (-1,2) {\textbf{B}};
			\draw[->] (0,-2) -- (0,1)node[left]{$f(u)$};
			\draw[->] (0,-2) -- (4,-2)node[below]{$u$};
			\draw[thick] (0,-2) .. controls (3,-2) and (1,0) .. (4,0);
			\draw[thick] (0,-2) -- (2,-2) node[below]{$u_c$};
			\draw[thick] (2,0) -- (4,0);
			\draw[thick,dotted] (2,-2) -- (2,0);
			\node[left] at (0,0) {$a$};
			\draw[->,thick] (1,-1)node[above]{\footnotesize original} -- (1.5,-1.5);
			\draw[->,thick] (1,0.5)node[above]{\footnotesize approximation} -- (2.5,0.1);
		\end{scope}
	\end{tikzpicture}

	\caption{Propagation of a front supported by a plane/line with a piecewise constant growth rate. \textbf{A}. The general geometry of the model. The 
	$\z$-axis can be absent. In this case, the model is considered on a plane and the growth happens on a line. \textbf{B}. Piecewise constant 
	approximation of a sigmoidal growth rate $f$ bound by the maximal growth rate $a$ at infinity.}

	\label{fig-problem}
\end{figure}
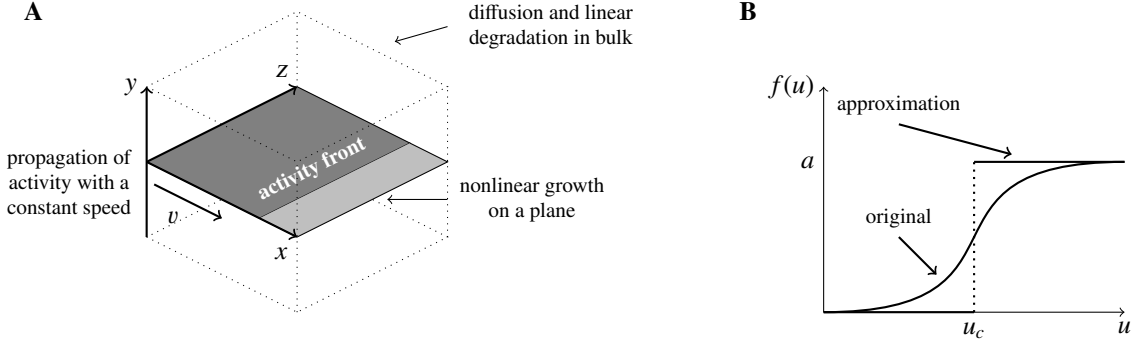

Let us consider a substance that is produced exclusively on the plane $\mathbb R^2(x,y)$ embedded in $\mathbb R^3(x,y,\z)$ with the production rate 
surface density $f(u)$ that depends on the local concentration $u$ of the substance. Equivalently, the line $\mathbb R(x)$ can be considered in $\mathbb 
R^2(x,y)$. Let us also assume that the substance can freely diffuse in the surrounding space with the diffusion coefficient $D$ and is degraded there 
with a linear rate constant $k$ (see Figure~\ref{fig-problem}A). Evolution of the concentration of the substance is then described by the following 
reaction-diffusion equation with a singular reaction term

\begin{equation}
\D_t u = D \Delta u -ku + f(u) \delta(y),
\label{eq3d}
\end{equation}

\noindent where $\Delta$ is the three- or two-dimensional Laplacian, $\delta$ is the Dirac $\delta$-function, and $f(u) \geqslant 0$ for $u \geqslant 
0$. Equation (\ref{eq3d}) is understood in the sense of generalised function (Schwartz's distributions), and $f$ and $u$ are assumed to be 
sufficiently well behaving for it to make sense. In particular, $u$ is assumed to be continuous.

Note that the initial value problem with equation (\ref{eq3d}) is equivalent to an initial-boundary problem for the linear equation $\D_t u = D \Delta 
u - ku$ in the half-space $y > 0$ with the nonlinear boundary condition $D \D_y u + f(u)/2 = 0$ on the boundary $y = 0$.

If the term $f(u)\delta(y)$ is replaced by the term $f(u)$, (\ref{eq3d}) becomes a classical reaction-diffusion equation with the reaction term $F(u) 
\deq f(u) - ku$. If $f(u)$ is monotone, increasing from 0 to some positive value at $[0,+\infty)$, and has a sigmoidal shape such that $F(0) = F(x_1) = 
F(x_2) = 0$ for some $x_1 > x_2 > 0$ and $F(x) < 0$ on $(0,x_2)$, while $F(x) > 0$ on $(x_2,x_1)$, then the corresponding nonsingular version of 
(\ref{eq3d}) belongs to the so called class of reaction-diffusion systems with a bistable nonlinearity (in the sense that if the Laplacian term is 
omitted, the resulting equation $\dot u = F(u)$ describes a bistable dynamical system). Such systems are known to produce planar travelling waves 
(fronts) of the form $u(t,x,\ldots) = \phi(x - \varv t)$, where $\ldots$ means spatial coordinates other than $x$ and $\varv$ is the constant 
travelling speed. A travelling solutions that connects the two stable steady states of the corresponding dynamical system, if exists, is unique (up to 
a spatial or temporal shift and a spatial rotation) with the unique value of the propagation velocity. Travelling waves in such bistable systems are 
well studied. Their utility is given by the fact that a sufficiently strong perturbation of the trivial state converges to such a solution at large 
time \cite{Kanel1962,Weinberger1975,Weinberger1978,Volpert1994}.

In some particular cases exact planar travelling solutions of the classical bistable reaction-diffusion equation can be obtained. One of such cases 
corresponds to the so called piecewise linear approximation of the reaction term. In this approach, a smooth reaction rate function $F$ is replaced by 
a piecewise linear (or rather affine) function that represents asymptotic regimes of $F$. For example, for $F(u) = f(u) - ku$ with a sigmoid $f$ such 
that $f'(0) = 0$ and $f(u) \to a$ at infinity a popular approximation corresponds to $F(u) = a\theta(u-u_c) - ku$, where $u_c$ is the single parameter 
of the approximation and $\theta$ is the Heaviside function. This piecewise linear functions corresponds to a piecewise constant approximation of the 
production term $f(u) = a\theta(u-u_c)$. A sketch of such an approximation is shown on Figure~\ref{fig-problem}B. The method of piecewise linear 
approximations to smooth functions in scientific modelling was pioneered by the Babylonian astronomer Naburimannu in about 500 BCE (from 610 BCE to 470 
BCE) \cite{vanderWaerden1974}. It was used for the right-hand sides of equations for dynamical systems as early as 1937 by the school of Andronov 
\cite{Andronov1981}. Such an approximation of the reaction term in models of propagation in reaction-diffusion systems was for the first time used in 
the work of McKean for a more complicated case of Nagumo's equation in neurophysics \cite{McKean1970}. It was later used for studying travelling waves 
in scalar bistable reaction-diffusion equations \cite{Wang2001,Petrovskii2005}, as well as in finding exact solutions for front propagation problems in 
models with more general dispersion kernels \cite{Nec2010,Volpert2010}. In most cases, piecewise approximations result in mathematically tractable 
exact solutions that have the same qualitative properties as general travelling solutions with regular reaction terms, at least where such properties 
are know.

This article considers equation (\ref{eq3d}) with a picewise constant approximation of $f(u)$ given by $a\theta(u-u_c)$. This substitution can be 
regarded as an approximation of a very thin transition zone of the sigmoid $f$, as was shown above. Thus, the chosen approximation is equivalent to a 
particular piecewise linear approximation of the reaction term in the classical reaction-diffusion problem. Such growth term may not even be 
approximate but exact if the production switches between two different states (``active'' and ``inactive'') as a function of the local concentration, 
as it is very commonly assumed in models of cell activity. The difference with the classical case is that now the degradation happens everywhere but 
the growth happens only on the special plane/line.

An analog of a travelling plane wave in the case of a surface-supported production, by the symmetry of the problem, is a concentration profile that 
spreads in some selected direction in the plane and that is translationary invariant in the perpendicular planar direction. If the $x$-axis is chosen 
along the propagation vector, the $y$-axis is chosen perpendicular to the supporting plane, then the solution of interest has the form $u(t,x,y,\z) = \varw(x - \varv t,y)$. In this situation, the problem becomes effectively two-dimensional. The equation of interest, thus, 
formally reads

\begin{equation}
\D_t u = D\Delta u - k u + a \theta(u - u^{}_c) \delta(y),
\label{main-equation}
\end{equation}

\noindent where $\Delta = \D_x^2 + \D_y^2$ is the two-dimensional Laplacian, and instead of the production supporting plane we have the production 
supporting line. Its exact meaning will be specified below.

This article presents an exact travelling along the $x$-axis solution of (\ref{main-equation}) of the form $u(t,x,y) = \varw(x - \varv t,y)$ with constant 
propagation velocity $\varv$ that connects the trivial stationary solution of the equation to a nontrivial one. Only the existence problem of an 
eternally propagating with constant speed solution and its uniqueness are studied (by an explicit construction). Questions of its stability and 
conditions on its initiation are not considered.

Let us reparametrise the variables of equation (\ref{main-equation}), so far only formally, to reduce the number of parameters. If we group together all the spatial 
independent variables $x^{}_i$ ($x_0 = x$, $x_1 = y$) such that they are rescaled simultaneously, the equations has four parameters and three variables 
for rescaling. In a generic case, this allows to reduce the number of parameters to one, using the appropriate reparametrisation. Assuming $a, k, D > 0$,with the mapping

\begin{equation}
t \mapsto t\frac{1}{k}, \quad\quad x^{}_i \mapsto x^{}_i \sqrt{\frac{D}{k}}, \quad\text{and}\quad u \mapsto u\frac{a}{\sqrt{kD}},
\label{rescaling}
\end{equation}

\noindent the equation takes the form

\begin{equation}
\D_t u = \Delta u - u + \theta(u - \alpha) \delta(y)
\label{eq-dimless}
\end{equation}

\noindent with a single parameter $\alpha = \frac{u^{}_c}{a} \sqrt{kD}$.

\section*{Results}

\begin{definition}
A real-valued continuous function $u \in \mathcal C(\Omega_t)\cap\mathcal S'(\Omega_t)$ is called a \emph{physical solution} of (\ref{eq-dimless}) if the level set $\{u|_{L_t} = \alpha\}$ has measure 0 and if $u$ respects in terms of $\mathcal S'$ the following equation

\begin{equation}
\D_t u  - \Delta u + u = \theta(u|_{L_t} - \alpha) \otimes \delta,
\label{eq-distributions}
\end{equation}

\noindent where the right-hand side is an element of $\mathcal S'(L_t) \otimes \mathcal S'(\mathbb R(y)) \subset \mathcal S'(\Omega_t)$.

\end{definition}

The definition is correct because $\theta$ is a regular distribution and can be seen as an element of $\mathcal L^1_\mathrm{loc}(\mathbb R)$. In addition 0 is the only point where the essential limit of $\theta$ is undefined. Thus, $\theta$ can be unambiguously composed with any measurable function $\phi\colon X \to \mathbb R$ such that $\phi^{-1}(0)$ has measure 0. The resulting object $\theta\circ \phi\in \mathcal L^1_\mathrm{loc}(X)$ is simply the class of the indicator function of $\phi^{-1}\big((0,+\infty)\big)$. Note that in this case the indicator of $\phi^{-1}\big([0,+\infty)\big)$ is a member of the same class.

The demand for $u$ to be tempered corresponds to a class of boundary conditions at infinity that exclude unphysically strong growth, which may be interpreted as the system in question not being strongly infuenced by the external world. Technically, this condition is enough to guarantee the uniqueness of the investigated travelling solution.

The main results of this article are represented by the three following theorems.

\begin{theorem}
\label{theorem1}

For any $\alpha \in (0,1/2)$ equation (\ref{eq-dimless}) has a unique nontrivial $x$-translation invariant stationary physical solution $u(t,x,y) = s(y)$ with

\begin{equation}
s(y) = \frac{e^{-|y|}}{2}.
\label{stationary}
\end{equation}

\end{theorem}

\begin{theorem}
\label{theorem2}

For any $\alpha \in (0,1/2)$ equation (\ref{eq-dimless}) has a travelling physical solution $u(t,x,y) = \varw(x - \varv t,y)$ with the front shape function

\begin{equation}
\varw(x,y) = \frac{1}{2\pi} \inth{0}{\infty}
K^{}_0\left(\sqrt{(x + \xi)_{}^2 + y_{}^2}\sqrt{1 + \frac{\varv_{}^2}{4}}\right) 
\exp\left(-\frac{\varv}{2}(x + \xi)\right)\, d\xi,
\label{u-v}
\end{equation}

\noindent where

\begin{equation}
\varv = 2 \ctg 2\pi \alpha
\label{v-dimless}
\end{equation}

\noindent is its propagation velocity.

\end{theorem}

\begin{theorem}
\label{theorem3}

The travelling physical solution from Theorem~\ref{theorem2} has the following properties:
\begin{enumerate}[label=(\roman*)]

\item For any $p \in \mathbb R$, $\varw(x,p)$ is a monotone decreasing function of $x$ on the whole $\mathbb R$, $\varw(p,y)$ is a monotone decreasing 
function of $y$ on $(0,+\infty)$ and a monotone increasing function on $(-\infty,0)$. \label{prop-monotone}

\item
The front $\varw$ connects the stationary states $0$ and $s$ in the following sense: for any $y$, $\lim\limits_{x\to 
+\infty}\varw(x,y) = 0$, $\lim\limits_{x\to -\infty}\varw(x,y) = s(y)$. \label{prop-connection}

\item
Travelling physical solution (\ref{u-v})--(\ref{v-dimless}) is unique (up to $x$-reflections and $x$-shifts of $\varw$) among travelling physical solutions with $\varw$ such that $\varw(x,0) < \alpha$ for all $x > x_1$ (respectively $x < x_1$) and $\varw(x,0) > \alpha$ for all $x < x_2$ (respectively $x > x_2$), where $x_1$ and $x_2$ are some real numbers.
\label{prop-unique}

\item
The propagation velocity $\varv$ can take any value from $(-\infty, +\infty)$. In particular, with $\alpha = 1/4$ the front is stationary. 
\label{prop-v}

\end{enumerate}
\end{theorem}

\subsection*{Proofs}

A preliminary lemma will be used in the proofs. Consider a differential operator over $\mathbb R^n(x_1,\ldots,x_n)$ with constant coefficients, \emph{viz.} $P(\D_{x_1},\ldots,\D_{x_n})$, where $P$ is a complex-valued polynomial of $n$ variables.

\begin{lemma}
\label{Polynomial}

If $P(-ix_1,\ldots,-ix_n) \neq 0$ for any $(x_1,\ldots,x_n) \in \mathbb R^n$, then $P(\D_{x_1},\ldots,\D_{x_n})$ has a unique tempered fundamental solution $\Phi$, \emph{viz.} a distribution $\Phi \in \mathcal S'$ that solves

\begin{equation*}
P(\D_{x_1},\ldots,\D_{x_n}) \Phi = \delta,
\end{equation*}

\noindent and for any $\Sigma \in \mathcal S'$ the linear PDE

\begin{equation}
P(\D_{x_1},\ldots,\D_{x_n}) u = \Sigma
\label{P-fundamental}
\end{equation}

\noindent has a unique tempered solution $u \in \mathcal S'$ given by

\begin{equation*}
u = \Phi * \Sigma.
\end{equation*}

\end{lemma}

\begin{proof}

Let us denote by $Q$ the complex-valued function on $\mathbb R^n$ given by $Q\colon (x_1,\ldots,x_n) \mapsto P(-ix_1,\ldots,-ix_n)$. Under the conditions of the lemma, (\ref{P-fundamental}) in $\mathcal S'$ is equivalent (by the Fourier transform and using the fact that the Fourier transform is an automorphism of $\mathcal S'$ as a linear space) to the following functional equation

\begin{equation}
Q \tilde u = \tilde \Sigma,
\label{P-functional}
\end{equation}

\noindent where $\tilde g \in \mathcal S'$ denotes the Fourier transform of $g \in \mathcal S'$.

Under the condition on $P$, $1/Q \in \mathcal O_M$ and thus $\tilde \Sigma/Q = (1/Q)\tilde \Sigma = \tilde \Sigma(1/Q)$ is well defined in $\mathcal S'$ \cite{Vladimirov}. It follows that $\tilde u = \tilde \Sigma/Q$ is the unique solution of (\ref{P-functional}) in $\mathcal S'$ seen as an equation on $\tilde u$. Indeed, let $\tilde u = \nu_1$ and $\tilde u = \nu_2$ solve (\ref{P-functional}). Consider any $\phi \in \mathcal S$. Then we have

\begin{equation}
0 = \langle Q(\nu_1 - \nu_2),\phi\rangle = \langle \nu_1 - \nu_2, Q\phi\rangle.
\end{equation}

\noindent But as $Q$ is nowhere vanishing, this is eqivalent to $\nu_1 = \nu_2$ in $\mathcal S'$. Thus the inverse Fourier transform of $\tilde \Sigma/Q$ is the unique solution to (\ref{P-fundamental}) in $\mathcal S'$.

Consider now the special case $\Sigma = \delta$ and thus $\tilde \Sigma = 1$. Then if we denote by $\Phi$ the inverse Fourier transform of $1/Q$, $\Phi$ is the unique fundamental solution of $P(\D_{x_1},\ldots,\D_{x_n})$ in $\mathcal S'$.

As $1/Q \in \mathcal O_M$, $\Phi * \Sigma$ exists in $\mathcal S'$ and is equal to the inverse Fourier transform of $\tilde \Sigma/Q$ \cite{Larcher2013}. This concludes the proof.

\end{proof}

\begin{proof}[Proof of Theorem~\ref{theorem1}]

A stationary $x$-translation invariant solution of (\ref{eq-dimless}) is a solution that depends only on $y$ (in terms of distributions it has the form $1(t,x)\otimes s(y)$, $s \in \mathcal S'(\mathbb R)$). Any such nontrivial phyisical solution corresponds to a continuous tempered solution of the following equation in $\mathcal S'(\mathbb R(y))$

\begin{equation}
\D^2_y s - s + \delta = 0.
\label{eq-s}
\end{equation}

In terms of Lemma~\ref{Polynomial}, we need to find a tempered fundamental solution of the differential operator $-\D^2_y + 1$, for which $Q(y) = y^2 + 1$ fulfills the requirement of the lemma. Therefore a known tempered solution to (\ref{eq-s}) given by $s(y) = e^{-|y|}/2$ is unique in $\mathcal S'$.

If $\alpha \geqslant 1/2 = s(0) = \max\limits_{y} s(y)$, though, $s$ is not a physical solution of (\ref{eq-dimless}).

\end{proof}

\begin{proof}[Proof of Theorem~\ref{theorem2}]

In the comoving frame finding a physical solution to equation (\ref{eq-dimless}) is equivalent to finding a real-valued continuous solution in $\mathcal S'(\Omega)$ with the level set $\{\varw|_L = \alpha\}$ of 0 measure in L of the following equation

\begin{equation}
-\Delta \varw -\varv \D_x \varw + \varw = \delta(y) \sum_{i} \chi^{}_{I_i}(x),
\label{general}
\end{equation}

\noindent where $\chi^{}_X$ is the indicator function of a subset $X \subset \mathbb R$ and $\{I_i\}$ is the set of nonintersecting closed intervals of 
$\mathbb R$ such that $\bigcup\limits_i I^{}_i = \{x \in \mathbb R \,|\, \varw(x,0) \geqslant \alpha\}$. The intervals $I^{}_i$ depend on $\varw$ rendering the equation nonlinear.

Consider problem (\ref{general}) as a nonhomogeneous generalised linear problem

\begin{equation}
- \Delta \varw - \varv \D_x \varw + \varw = \Sigma,
\label{general2}
\end{equation}

\noindent with a given source term $\Sigma \in \mathcal S'$. The polynomial $x^2 + y^2 + i\varv x + 1$ does not vanish in $\mathbb R^2$, threefore by Lemma~\ref{Polynomial} there is a unique tempered solution of (\ref{general2}) given by $\Phi*\Sigma$, where $\Phi$ is the unique tempered fundamental solution of the operator $A = -\Delta -\varv \D_x + 1$. By introduction of a new dependent variable $\Phi = \Psi \exp(-\varv x /2)$ and taking into account the identity $\exp(\eta x)\delta(x,y) = \delta(x,y)$, $\eta \in \mathbb R$, the problem $A\Phi = \delta$ is transformed to

\begin{equation*}
-\Delta \Psi + \left(1 + \frac{\varv_{}^2}{4}\right) \Psi = \delta.
\end{equation*}

A tempered fundamental solution of the operator $-\Delta + (1 + \varv^2/4)$ in this equation, which is unique by Lemma~\ref{Polynomial}, is known to be

\begin{equation*}
\Psi(x,y) = \frac{1}{2\pi} K^{}_0\left(\sqrt{x_{}^2 + y_{}^2}\sqrt{1 + \frac{\varv_{}^2}{4}}\right),
\end{equation*}

\noindent where $K^{}_0$ is the Macdonald function. Therefore,

\begin{equation*}
\Phi(x,y) = \frac{1}{2\pi} K^{}_0\left(\sqrt{x_{}^2 + y_{}^2}\sqrt{1 + \frac{\varv_{}^2}{4}}\right) \exp\left(-\frac{\varv}{2}x\right),
\label{fundamental-phi}
\end{equation*}

\noindent is a fundamental solution of $A$. It is easy to check, using the asymptotics of the Macdonald function $K^{}_0(x) \sim 
\dfrac{e^{-x}}{\sqrt{x}}$ at infinity, that $\Phi$ is tempered and thus unique tempered solution by Lemma~\ref{Polynomial}.

Thus the unique tempered solution of (\ref{general2}) $\varw = \Phi * \Sigma$ provides for (\ref{general}) a necessary condition on a travelling solution

\begin{equation}
\varw(x,y) = \frac{1}{2\pi} \inth{\bigcup\limits_i I_i}{}
K^{}_0\left(\sqrt{(x - \xi)_{}^2 + y_{}^2}\sqrt{1 + \frac{\varv_{}^2}{4}}\right) 
\exp\frac{\varv(\xi - x)}{2}\, d\xi,
\label{u-expression}
\end{equation}

\noindent where again the set $\{I^{}_i\}$ depends on $\varw$.

Let us additionally impose the following condition on the solution to be found: $\varw(x,0) > \alpha$ for all $x < x_1$, $\varw(x,0) < \alpha$ for all $x > x_2$, for some $x_1, x_2 \in \mathbb R$, which is motivated by the analogy with the calssical bistable travelling wave connecting the trivial state with a nontrivial state of the well mixed system.

This condition translates in $\bigcup\limits_i I_i \neq \varnothing$, $\sup \bigcup\limits_i I_i < +\infty$, and in the existence of an infinitely large interval $J = (-\infty, -\beta] \in \{I_i\}$. Without a loss of generality we can assume that $\sup \bigcup\limits_i I_i = 0$ and thus $\beta \geqslant 0$.

Suppose that $\bigcup\limits_i I^{}_i \neq (-\infty,0]$. It means that $\beta > 0$ and $J \neq (-\infty,0]$. By definition of $I^{}_i$ and by continuity of $\varw$ we must have $\varw(-\beta,0) = \varw(0,0)$, and thus the following equality

\begin{equation*}
\intL{-\infty}{-\beta}
K^{}_0\big(-b(\xi+\beta)\big)\exp\left(\frac{\varv}{2}(\xi+\beta)\right)\,d\xi\\
+
\intL{(\bigcup\limits_i I^{}_i)\setminus J}{}
K^{}_0\big(b(\xi+\beta)\big)\exp\left(\frac{\varv}{2}(\xi+\beta)\right)\,d\xi\\
=
\intL{\bigcup\limits_i I^{}_i}{}
K^{}_0(-b\xi)\exp\frac{\varv\xi}{2}\,d\xi,
\end{equation*}

\noindent Where, for brevity, $b = \sqrt{1 + \varv^2/4}$. By the premise, $\bigcup\limits_i I_i \subset (-\infty,0]$ but $\bigcup\limits_i I_i \neq 
(-\infty,0]$, and both these sets are closed. This implies

\begin{equation*}
\intL{-\infty}{0}
K^{}_0(-b\xi)\exp\frac{\varv\xi}{2}\,d\xi
>
\intL{\bigcup\limits_i I^{}_i}{}
K^{}_0(-b\xi)\exp\frac{\varv\xi}{2}\,d\xi.
\end{equation*}

\noindent Because all integrands are positive, this contradicts the previous relation. Therefore, the premise is false, and $\bigcup\limits_i I_i = J = (-\infty,0]$.

As a consequence, (\ref{u-expression}) is greatly simplified. The profile of the travelling solution is given by

\begin{equation}
\varw(x,y) = \frac{1}{2\pi} \inth{0}{\infty}
K^{}_0\left(\sqrt{(x + \xi)_{}^2 + y_{}^2}\sqrt{1 + \frac{\varv_{}^2}{4}}\right) 
\exp\left(-\frac{\varv}{2}(x + \xi)\right)\, d\xi,
\end{equation}

\noindent where the velocity is implicitly defined by $\varw(0,0) = \alpha$

\begin{equation}
2\pi \alpha = \inth{0}{\infty}\exp\left(-\frac{\varv \xi}{2}\right)K^{}_0\left(\xi \sqrt{1 + \frac{\varv^2_{}}{4}}\right)\,d\xi.
\label{v-implicit}
\end{equation}

The integral on the right-hand side of this expression can be computed in a close form giving an explicit formula for the propagation velocity. Indeed,
using the table integral

\begin{equation*}
\intl{0}{\infty} e^{-c x} K^{}_0(b x)\,dx = \frac{\arccos \frac{c}{b} }{\sqrt{b_{}^2 - c_{}^2}},
\end{equation*}

\noindent where in our case $c = \varv/2$ and $b = \sqrt{1 + \varv_{}^2/4}$, we can transform (\ref{v-implicit}) to

$$
2\pi \alpha = \arccos \frac{1}{\sqrt{\frac{4}{\varv_{}^2} + 1}},
$$

\noindent from where, after a simple rearrangement, we obtain

\begin{equation}
\varv = 2 \ctg 2\pi \alpha.
\end{equation}

\end{proof}

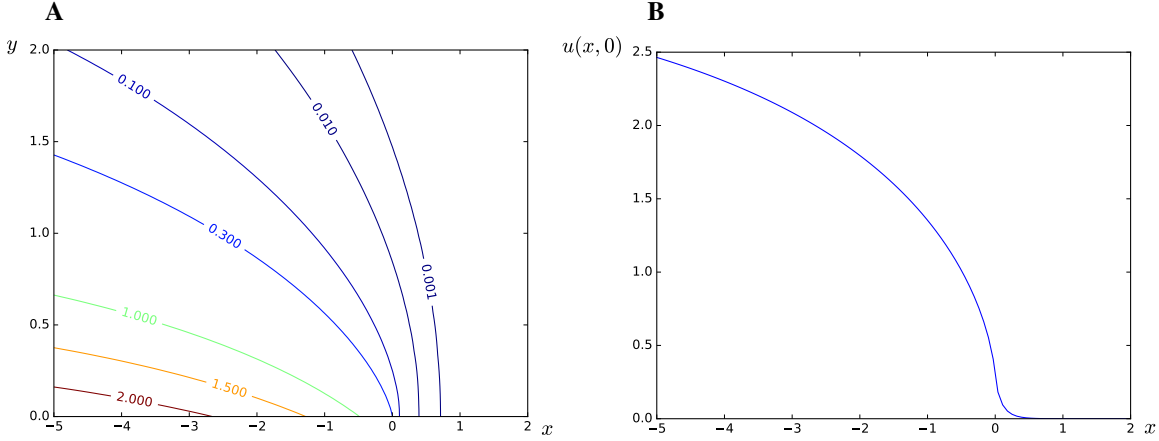
\begin{figure}[t]
\begin{tikzpicture}
	\node at (0,0) {
		\includegraphics[scale=0.4]{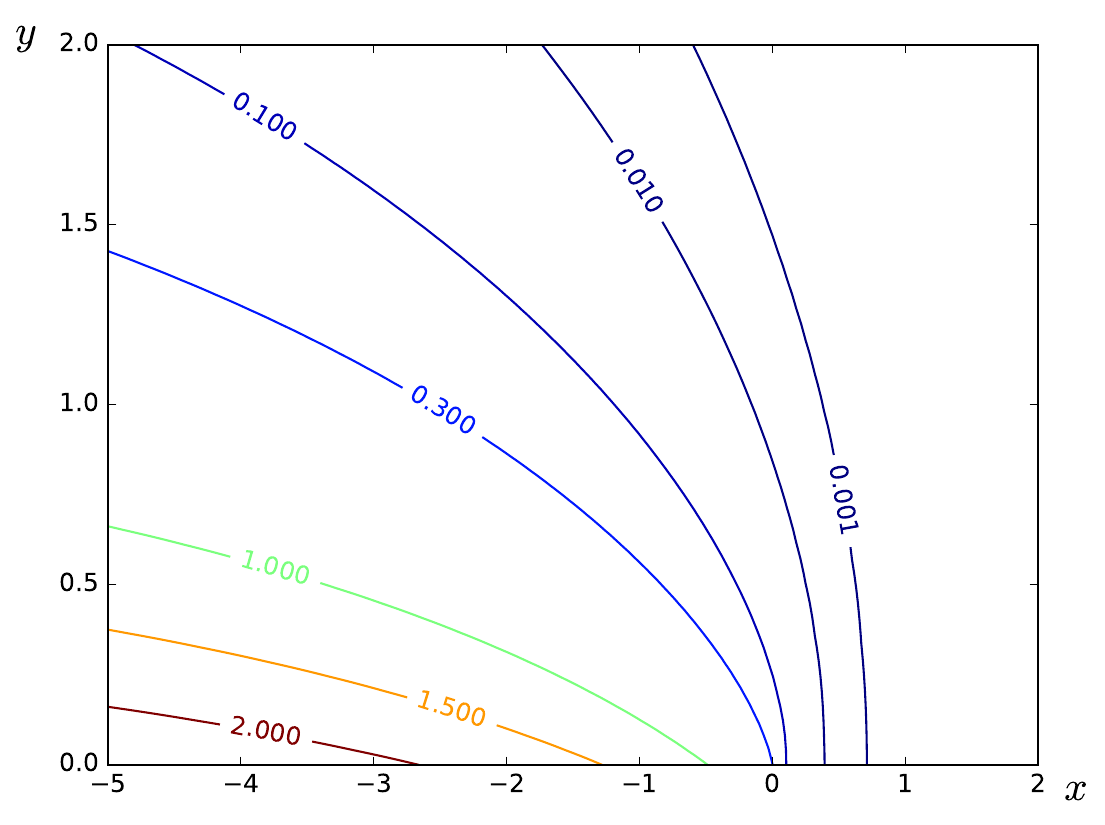}
	};
	\node at (-3,3) {\bf A};
	\node at (7.7,0) {
		\includegraphics[scale=0.4]{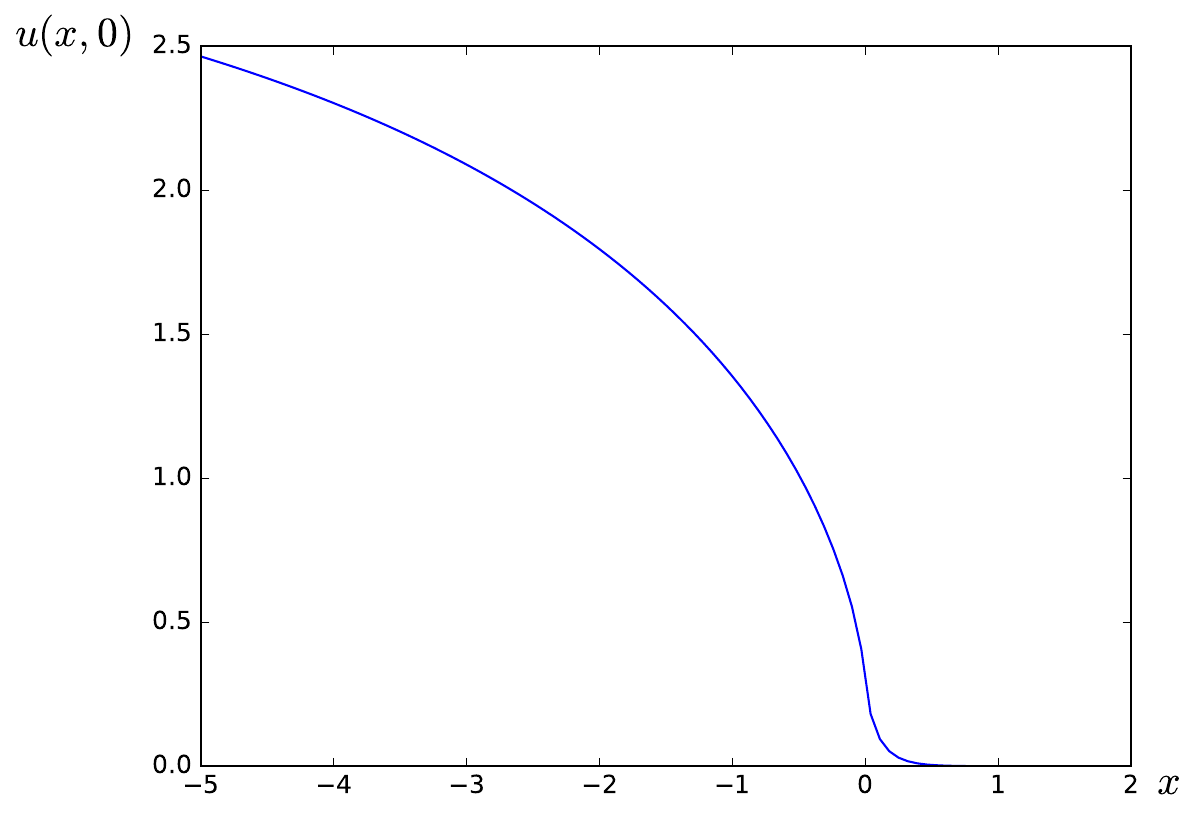}
	};
	\node at (5,3) {\bf B};
	\fill[color=white] (3.5,2.3) rectangle (4.6,2.8);
	\node at (4.1,2.57) {$\varw(x,0)$};
\end{tikzpicture}
\caption{
	The front $\varw$ of the travelling solution $u(x,y,t) = \varw(x-\varv t,y)$ of (\ref{main-equation}) given by (\ref{sol-units}) for $a = 2\pi$, $k = D = 1$, and $u_c = 0.3$ ($\varv \approx 6.47$). 
	\textbf{A}. A contour plot of $\varw$ in the $xy$-plane ($y \geqslant 0$). Note that $\varw(x,y) = \varw(x,-y)$. \textbf{B}. The value of $\varw$ on the $x$-axis. 
	Here the production is active on the $x$-axis with $x \in (-\infty,0]$ and inactive everywhere else.
}
\label{fig-profile}
\end{figure}

To recover the expression of the front propagation velocity in the original units it is enough, as (\ref{rescaling}) suggests, to multiply this 
expression by $\sqrt{kD}$. So, in the original units we have

\begin{equation}
\varv = 2\sqrt{kD} \ctg \frac{2\pi u^{}_c \sqrt{kD}}{a}
\label{v-units}
\end{equation}

and the expression for the wave profile becomes (see Figure~\ref{fig-profile})

\begin{equation}
\varw(x,y) = \frac{a}{2\pi D} \inth{0}{\infty}
K^{}_0\left(\sqrt{(x + \xi)_{}^2 + y_{}^2}\frac{\sqrt{4kD + \varv_{}^2}}{2D}\right) 
\exp\left(-\frac{\varv}{2D}(x + \xi)\right)\, d\xi.
\label{sol-units}
\end{equation}

\begin{proof}[Proof of Theorem~\ref{theorem3}]

The part of property \ref{prop-monotone} with $x$ is obvious once (\ref{u-v}) is rewritten in the form

\begin{equation}
\varw(x,y) = \frac{1}{2\pi} \inth{x}{+\infty}
K^{}_0\left(\sqrt{\xi_{}^2 + y_{}^2}\sqrt{1 + \frac{\varv_{}^2}{4}}\right) 
\exp\left(-\frac{\varv}{2}\xi\right)\, d\xi.
\label{alternative}
\end{equation}

The part with $y$ follows from the monotonicity of $K^{}_0(x)$ for $x \in (0,+\infty)$.

For property \ref{prop-connection}, the limit at $x \to +\infty$ is clear from (\ref{alternative}). If we again denote for brevity $b = \sqrt{1 + 
(\varv/2)_{}^2}$ and $c = \varv/2$ and use one of the integral representations of the Macdonald function, the limit at $x \to -\infty$ is equal to

\begin{multline}
\lim_{x\to -\infty} \varw(x,y) = \frac{1}{2\pi} \inth{-\infty}{+\infty}
K^{}_0\left(b\sqrt{\xi_{}^2 + y_{}^2}\right) 
e_{}^{-c\xi}\, d\xi =
\frac{1}{4\pi} \inth{-\infty}{+\infty} \inth{0}{+\infty}
\frac{1}{t}\exp\left(-c\xi-t-\frac{b^2\xi^2}{4t}-\frac{b^2y^2}{4t}\right)\, dt d\xi = \\
= \frac{1}{4\pi}\inth{0}{+\infty}\frac{1}{t}\exp\left(-t-\frac{b^2y^2}{4t}\right)
\inth{-\infty}{+\infty}\exp\left(-c\xi-\frac{b^2\xi^2}{4t}\right)\,d\xi dt =
\frac{1}{2b\sqrt{\pi}}\inth{0}{+\infty}\frac{1}{\sqrt{t}}\exp\left(-\frac{t}{b^2}-\frac{b^2y^2}{4t}\right)\,dt = \\
= \frac{1}{2\sqrt{\pi}}\inth{0}{+\infty}\frac{1}{\sqrt{t}}\exp\left(-t-\frac{y^2}{4t}\right)\,dt = \frac{1}{2}e^{-|y|} = s(y).
\label{connection}
\end{multline}

Property \ref{prop-unique} follows from the explicit construction and Lemma~\ref{Polynomial}.

Property \ref{prop-v} is obvious from (\ref{v-dimless}).

\end{proof}

The value of the propagation velocity given by (\ref{v-dimless}) is a periodic function of $\alpha$ with the period of $1/2$, which monotonely 
decreases from $+\infty$ to $-\infty$ on $(0,1/2)$. One may wonder what happens if $\alpha$ surpasses this interval. As it is seen from 
(\ref{connection}), we have $\sup s = 1/2$. As the travelling solution is monotone in $x$, we must have $\sup \varw = 1/2$. Therefore, if $\alpha
> 1/2$, the source term in (\ref{eq-dimless}) is equal to 0. This means, that there is neither the nontrivial $x$-translation invariant stationary
solution nor the travelling one, if $\alpha$ is too high.

These properties well parallel the properties of the travelling solutions in classical reaction-diffusion systems with bistable nonlinearities 
\cite{Weinberger1975,Volpert1994}. There is, however, no simple geometric phase space analysis available to support them (like the one shown on 
Figure~\ref{fig-sketches}A, see \cite{Wang2001} for details). Another difference is in the asymptotics of the travelling wave in its leading front. The 
classical case has an exponentially decaying front, while for solution (\ref{u-v}), $\varw(x,0)$ decays as $e^{-\gamma x}/\sqrt{x}$ at $x \to +\infty$, 
where $\gamma$ is a positive constant (this follows from (\ref{alternative}) and the asymptotic properties of $K_0$).

\section*{Discussion and Conclusion}

In addition to the properties from Theorem~\ref{theorem3}, which parallel the ones of the classical bistable reaction-diffusion equation, solution (\ref{v-units})--(\ref{sol-units}) has a property specific to this particular problem. Expression (\ref{v-units}) for the propagation speed has a well defined limit at $k \to 0$: $\varv_0 = a/(\pi u_c)$. The corresponding limit of (\ref{sol-units}) diverges at any $y$ in the 
limit $x \to -\infty$ and thus connects the trivial state with infinity. This limit profile is in fact a travelling solution of 
(\ref{main-equation}) with no degradation. Indeed, with $k = 0$ and with a reparametrisation different from (\ref{rescaling}), \emph{viz.} $u \mapsto 
au/D$ and $t \mapsto t/D$, the equation becomes $\D_t u = \Delta u + \theta(u - \alpha)\delta(y)$, where now $\alpha = D u_c/a$. Following the same 
route, it is possible to show that all the subsequent arguments are valid with the replacement of the expression $1 + \varv^2/4$ by $\varv^2/4$. If 
$\varv \leqslant 0$, the integral in condition (\ref{v-implicit}), that now reads as $2\pi\alpha = 
\intl{0}{\infty}e^{-\varv\xi/2}K_0(|\varv|\xi/2)\,d\xi$, diverges. However, under assumption $\varv > 0$, this condition gives the value $\varv = 
(\pi\alpha)^{-1}$, which becomes $\varv_0$ in the original units, as expected. Furthermore, the corresponding solution $u$ in the original units 
exactly equals to the limit of (\ref{sol-units}) at $k \to 0$. Note however that the premise of Lemma~\ref{Polynomial} fails in this case and the uniqueness of a monotone travelling solution cannot be guaranteed.

\begin{figure}[!t]
\centering
\begin{tikzpicture}
	\begin{scope}
		\node at (0,2) {\bf A};
		\draw[->,line width=0.5pt] (0,-1) -- (0,1)node[above left]{$\D_x \varw$};
		\draw[->,line width=0.5pt] (-1,0) -- (2,0)node[below right]{$\varw$};
		\draw[dotted,line width=0.5pt] (0.5,1) -- (0.5,-1)node[below]{$\varw = u_c$};
		\node at (0,0) {$\bullet$};
		\node at (1.75,0) {$\bullet$};
		\draw[->-=0.5,very thick] (1.75,0) -- (0.5,-0.5);
		\draw[->-=0.5,line width=0.5pt] (2,0.25*0.5/1.25) -- (1.75,0);
		\draw[->-=0.5,line width=0.5pt] (1.75,0) -- (2,-0.25);
		\draw[->-=0.5,line width=0.5pt] (0.6,1.25*5/6) -- (1.75,0);
		\draw[->-=0.5,very thick] (0.5,-0.5) -- (0,0);
		\draw[->-=0.5,line width=0.5pt] (-1,1) -- (0,0);
		\draw[->-=0.5,line width=0.5pt] (0,0) -- (0.5,0.5*0.5/1.25);
		\draw[->-=0.5,line width=0.5pt] (0,0) -- (-1,-1*0.5/1.25);
	\end{scope}

	\begin{scope}[xshift=4cm]
		\node at (0,2) {\bf B};
		\draw[->,line width=0.5pt] (0,-1) -- (0,1)node[above left]{$\D_x \varw$};
		\draw[->,line width=0.5pt] (-1,0) -- (2,0)node[below right]{$\varw$};
		\draw[dotted,line width=0.5pt] (1,1) -- (1,-1)node[below]{$\varw = u_c$};
		\draw[very thick,dotted] (-1,0) -- (1,0);
		\begin{scope}
			\clip (-1,-1) rectangle (1,1);
			\foreach \x in {-1.2,-0.6,...,3} {
				\draw[->-=0.75-0.2*\x,line width=0.5pt] (\x-2,1) -- (\x,0);
				\draw[->-=0.23+0.28*\x,line width=0.5pt] (\x+0.2,-1) -- (\x-1.8,0);
			}
			\node at (0,0) {$\bullet$};
			\draw[very thick,->-=0.75] (2,-1) -- (0,0);
		\end{scope}
		\draw[very thick,->-=0.5] (2,-0.5) -- (1,-0.5);
		\draw[->-=0.5,line width=0.5pt] (2,-1) to[in=-5,out=160] (1,-0.8);
		\draw[->-=0.5,line width=0.5pt] (1,0.4) to[in=90,out=-30] (1.4,0) to[in=5,out=-90] (1,-0.2);
		\draw[->-=0.5,line width=0.5pt] (1,0.7) to[in=90,out=-30] (1.7,0) to[in=2,out=-90] (1,-0.3);
	\end{scope}

	\begin{scope}[xshift=8cm]
		\node at (0,2) {\bf C};
		\draw[->,line width=0.5pt] (0,-1) -- (0,1)node[above left]{$\D_x \varw$};
		\draw[->,line width=0.5pt] (-1,0) -- (2,0)node[below right]{$\varw$};
		\draw[dotted,line width=0.5pt] (0.5,1) -- (0.5,-1)node[below]{$\varw = u_c$};
		\node at (0,0) {$\bullet$};
		\node at (1.75,0) {$\bullet$};
		\draw[->-=0.5,line width=0.5pt] (1.75,0) -- (0.5,-0.75);
		\draw[->-=0.5,line width=0.5pt] (2,0.25*0.75/1.25) -- (1.75,0);
		\draw[->-=0.5,line width=0.5pt] (1.75,0) -- (2,-0.25*0.75/1.25);
		\draw[->-=0.5,line width=0.5pt] (0.5,0.75) -- (1.75,0);
		\draw[->-=0.5,very thick] (0.5,-0.5*0.75/1.25) -- (0,0);
		\draw[->-=0.5,very thick] (0,0) -- (0.5,0.5*0.75/1.25);
		\draw[->-=0.5,line width=0.5pt] (-1,0.75*0.75/1.25) -- (0,0);
		\draw[->-=0.5,line width=0.5pt] (0,0) -- (-1,-0.75*0.75/1.25);
		\draw[->-=0.5,very thick] (0.5,0.5*0.75/1.25) to[out=-10,in=90] (0.8,0) to[out=-90,in=10] (0.5,-0.5*0.75/1.25);
	\end{scope}

\end{tikzpicture}

\caption{Sketches of travelling solution profiles $\varw$ of the classical piecewise linear equation in the original units on the phase plane $(\varw,\D_x \varw)$. The 
	trajectory that corresponds to a front is depicted as the thick line. \textbf{A}. The regular travelling front that connects the trivial and the 
	nontrivial steady states for the case $\varv > 0$. It is the standard heteroclinic trajectory on the phase plane. Only separatrices of the saddles 
	are shown. See \cite{Wang2001} for the details. \textbf{B}. The travelling solution with $k = 0$ that connects the trivial steady state and infinity with bound derivative at infinity. The case is degenerated and the whole $\varw$-axis consists of steady states for 
	$\varw < u_c$. A generic smooth case would correspond to a saddle-node bifurcation at the origin in this situation. \textbf{C}. The homoclinic 
	stationary solution that connects the trivial steady state with itself in the case when the monotone front 
	travels with $\varv > 0$.
}
\label{fig-sketches}
\end{figure}
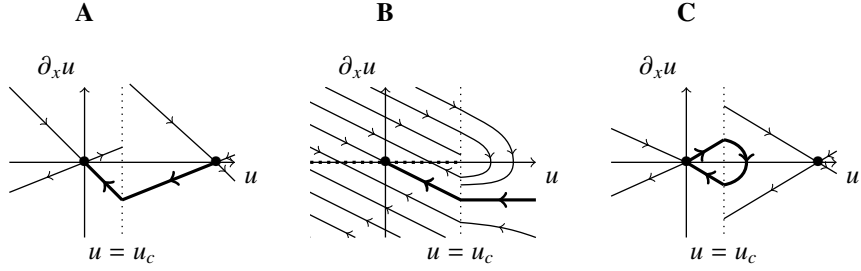

The existence of a travelling solution with $k = 0$ is not the unique property of the codimension-1 supported reaction-diffusion system. The specific 
property is that the travelling velocity $\varv_0$ does not depend on the intensity of diffusion ($D$), unlike the corresponding velocity for the 
classical case. The classical case is described by the equation $\D_t u = D\D_x^2 u + a\theta(u - u_c)$ and results in a similar travelling solution 
(see Figure~\ref{fig-sketches}B), which connects the trivial homogeneous solution with infinity and which is characterized by the velocity 
$\sqrt{aD/u_c}$ (this can be seen as the limit of the unique travelling solution for $k \to 0$ from \cite{Petrovskii2005}).

One should be careful about the last point. The travelling solution that connects the trivial solution with infinity in both classical and 
codimension-1 supported cases with piecewise constant production is unique (if we demand bound derivatives at infinity for physical relevance). This 
is, however, a special property of the piecewise linear approximation. At least in the classical case with a generic smooth monotone growth function 
$f$ such that $f(0) = f'(0) = 0$, $f''(0) > 0$, and $f(u) \to \mathrm{const}$ at $u \to +\infty$, one can show that this solution corresponds to the 
solution with the minimal propagation velocity. This situation is typical for FKPP-type (monostable) reaction terms \cite{Volpert1994}, to which the 
bistable case degenerates with $k \to 0$ locally near $u = 0$. From what is known for the studied FKPP cases, we can expect that this slowest solution 
corresponds to the long term limiting regime of generic perturbation of the trivial solution. On the other hand, the independence of $\varv_0$ of $D$ 
is a generic property that does not depend on the piecewise linear approximation. Indeed, it simply follows from the scaling properties of 
(\ref{eq3d}).

On one hand, the significance of the exact solution presented in the current article is limited. Indeed, it applies only to a very specially 
degenerated case. The approach to solve this case is not generalisable to generic production terms $f$ in (\ref{eq3d}). It relies on the reduction of 
the problem to a linear problem with generalised sources and the uniqueness of the travelling solution follows from the condition (\ref{v-implicit}). 
On the other hand, reasonably assuming that it is a limit case of some family of general regular cases, it hints on the expected properties in a 
generic case.

Having in mind this point, it is worth to note that for any $\alpha < 1/4$ (and thus for the advancing front with $\varv > 0$), equation 
(\ref{eq-dimless}) has a stationary solution in addition to the travelling one. Its existence can be established with the same methods from 
(\ref{u-v}). It is supported by the production on a finite interval, so $\varw(x,0) \geqslant \alpha$ for $x$ on some interval $[-\zeta,\zeta]$ (up to translations), 
where $\zeta$ depends on $\alpha$. This solution corresponds to the homoclinic stationary solution in the classical reaction-diffusion case (see 
Figure~\ref{fig-sketches}C).

The regularity of solution (\ref{u-v})--(\ref{v-dimless}), as expected, worsens at the special line and at the point of transition from no activity on the line to production. Being smooth (and even analytic) in $\Omega_t\setminus L_t$ it is only Lipschitz-continuous at $L_t\setminus\{x=\varv t, y = 0\}$ and only H\"older-continuous (with any exponent in $(0,1)$) at $\{x=\varv t,y=0\}$. The former loss of regularity will be observed for more regular reaction terms $f$, while the latter is specific only for the step-function production.

Interestingly, an attempt to use the solution method outlined in this article for the case of a reaction supported by a line in space, that is when in 
(\ref{main-equation}) the Laplacian is three-dimensional and the production term is replaced by $a\theta(u-u_c)\delta(y,\z)$, fails. The problem is in 
the divergence of any integral analogous to (\ref{v-implicit}). The divergence is accumulated at an infinitesimal distance. The same problem will be 
encountered with a general model of (\ref{eq3d}), where the production term is $f(u)\delta(y,\z)$ and the value $f(u)$ and the meaning of the whole 
term are not defined, as $u$ must be divergent at the production line. This either means that the approximation of an infinitely thin production zone 
is incorrect, that a different formalization is required to make sense of a function of a singular distribution (like, perhaps, the Colombeau 
formalism), or that there is no travelling solution with a constant velocity in this case.

All in all, there seem to be some deep structure in various cases of bistable reaction-diffusion systems when it comes to travelling solutions that is 
reflected by their common properties, which is also suggested by the work of Fang and Zhao on very general monotone semiflows with a bistable 
structure \cite{Fang2015}.

\section*{Acknowledgements}

The author is grateful to Danielle Hilhorst for a critical discussion of the mathematical content of the work and to Fyodor V. Tkachev for pointing out 
the piecewise linear approximations in Babylonian astronomy.

\blfootnote{\emph{e-mail:} \texttt{anton.zadorin@mis.mpg.de}}

\end{document}